\documentclass[12pt, a4paper]{amsart}
\usepackage[hmargin=2.5cm, vmargin=2cm, includefoot, twoside]{geometry}
\usepackage{latexsym, amssymb, mathrsfs, graphics}
\usepackage[latin1]{inputenc}  
\usepackage[T1]{fontenc}       

\theoremstyle{plain}
\newtheorem{theointro}{Theorem}

\newtheorem{corintro}[theointro]{Corollary}

\newtheorem{thm}{Theorem}[section]
\newtheorem{cor}[thm]{Corollary}
\newtheorem{prop}[thm]{Proposition}
\newtheorem{lem}[thm]{Lemma}
\newtheorem{claim}{Claim}
\newtheorem*{claim*}{Claim}
\newtheorem*{thm*}{Theorem}
\newtheorem*{prop*}{Proposition}

\theoremstyle{definition}

\newtheorem{rem}[thm]{Remarks}

\newcommand{\N}{\mathbb{N}}
\newcommand{\Z}{\mathbb{Z}}
\newcommand{\la}{\langle}
\newcommand{\ra}{\rangle}
\newcommand{\inv}{^{-1}}
\newcommand{\MF}{\mathscr{M}(F)}
\newcommand{\Meucl}{\mathscr{M}_{\mathrm{Eucl}}(F)}
\DeclareMathOperator{\proj}{proj}

\DeclareMathOperator{\algrk}{alg-rk}\DeclareMathOperator{\rk}{rk}
\DeclareMathOperator{\Isom}{Isom}

\newcommand{\field}[1]{\mathbb{#1}}
\newcommand{\integers}{\ensuremath{\field{Z}}}

\newcommand{\naturals}{\ensuremath{\field{N}}}

\begin{document}


\title[On geometric flats in {${\rm CAT}(0)$} Tits buildings]
{On geometric flats in the {${\rm CAT}(0)$} realization of Coxeter groups and Tits buildings}%

\author[P.-E. Caprace]{Pierre\,-Emmanuel Caprace$^*$}
\thanks{$^*$ F.N.R.S. research fellow}
\address{D\'epartement de Math\'ematiques\\
Universit\'e libre de Bruxelles, CP216\\
Bd du Triomphe, 1050 Bruxelles (Belgium).}
\email{pcaprace@ulb.ac.be}

\author[F. Haglund]{Fr\'ed\'eric Haglund}
\address{Laboratoire de math\'ematiques\\
Universit\'e Paris-Sud - B\^at. 425\\
91405 Orsay Cedex (France).}
\email{frederic.haglund@math.u-psud.fr}



\date{\today}
\begin{abstract}
Given a complete CAT(0) space $X$ endowed with a geometric action of a group $\Gamma$, it is known that if
$\Gamma$ contains a free abelian group of rank $n$, then $X$ contains a geometric flat of dimension $n$. We
prove a converse of this statement in the special case where $X$ is a convex subcomplex of the CAT(0)
realization of a Coxeter group $W$, and $\Gamma$ is a subgroup of $W$. In particular a convex cocompact subgroup
of a Coxeter group is Gromov-hyperbolic if and only if it does not contain a free abelian group of rank 2. Our
result also provides an explicit control on geometric flats in the CAT(0) realization of arbitrary Tits
buildings.
\end{abstract}

\maketitle

%




\section*{Introduction}

Let $X$ be a complete ${\rm CAT}(0)$ space and $\Gamma$ be a group
acting properly discontinuously and cocompactly on $X$. It is a
well known consequence of the so called flat torus theorem (see
\cite[Corollary II.7.2]{BH99}) that:

\smallskip  \emph{$({\Z}^n\Rightarrow {\field E}^n)$: if $\Gamma$ contains a free abelian group of rank
$n$, then $X$ contains a geometric flat of dimension $n$.}

\smallskip
Recall that a \textbf{(geometric) flat} of dimension $n$, also called \textbf{(geometric) $n$--flat},  is a
closed convex subset of $X$ which is isometric to the Euclidean $n$-space. One may wonder whether a converse of
this statement does hold, that is to say, whether the presence of a geometric $n$--flat in $X$ is reflected in
$\Gamma$ by the existence of a free abelian group of rank $n$. This question goes back at least to Gromov
\cite[\S$6.B_3$]{Gro93}.



 In the case $n=2$, in view of the flat plane theorem
(see \cite[Corollary III.H.1.5]{BH99}), this question can be
stated as follows:

\smallskip  \emph{If $X$ is not hyperbolic, does $\Gamma$
contains a copy of $\Z \times \Z$?}

\smallskip  The answer is known to be positive in the
following cases:
\begin{itemize}
\item  $\Gamma$ is the fundamental group of a closed aspherical
$3$-manifold, see \cite{KK04}.

\item $X$ is a square complex satisfying certain technical
conditions, see \cite{Wi05}.
\end{itemize}

A  \textbf{combinatorially convex subcomplex} of the Davis complex  $|W|_0$  of a Coxeter group $W$ is an
intersection of closed half-spaces of  $|W|_0$. The following result shows that, if $X$ is a such a
combinatorially convex subcomplex of
 $|W|_0$,  and if $\Gamma\subset W$ acts
cellularly, then the converse of the property $({\field Z}^n\Rightarrow {\field E}^n)$ above holds
for all $n$:

\begin{theointro}\label{thmA}
Let $X$ be a combinatorially convex subcomplex of the Davis complex $|W|_0$ of a Coxeter group $W$.  Let $\Gamma$ be a subgroup
of $W$ which preserves $X$ and whose induced action on $X$ is cocompact. If $X$ contains a geometric $n$--flat,
then $\Gamma$ contains a free abelian group of rank $n$.
\end{theointro}

Since half-spaces are CAT(0)-convex, combinatorially convex subcomplexes are CAT(0)-convex as well. We do not know if the theorem above is still true when  $X$ is only assumed to be a CAT(0)-convex subset of $|W|_0$. We note that in general the intersection $\bar X$ of the closed half-spaces of $|W|_0$ containing $X$ is not cocompact under $\Gamma$. Yet $\Gamma$ is still cofinite on the set of walls separating $X$ (or $\bar X$), and perhaps this is enough.

\begin{corintro}\label{thmA'}
Let $X$ be a CAT(0) convex subcomplex of the Davis complex $|W|_0$ of a Coxeter group $W$.  Let $\Gamma$ be a subgroup
of $W$ which preserves $X$ and whose induced action on $X$ is cocompact. If $X$ contains a geometric $n$--flat,
then $\Gamma$ contains a free abelian group of rank $n$.

\end{corintro}
\begin{proof}

The Corollary follows by Theorem~\ref{thmA} because, since $X$ is a subcomplex, the intersection of the closed
half-spaces of $|W|_0$ containing $X$ is a combinatorially convex $\Gamma$-cocompact subcomplex $\bar X$.

We sketch  the argument. The key-point is that $X^0$ is convex for the combinatorial distance. First, any two
vertices $x,y$ of $X$ may be joined by a combinatorial geodesic $(x_0=x,\dots,x_n=y)$ all of whose vertices
belong to the smallest subcomplex of $|W|_0$ containing the CAT(0) geodesic between $x$ ad $y$ (see
\cite[Lemme~4.9]{HP98}). Since $X$ is a CAT(0) convex subcomplex, it follows that $x_0,\dots,x_n$ belong to
$X^0$. Now any combinatorial geodesic between $x,y$ may be joined to $(x_0,\dots,x_n)$ by a sequence of
geodesics, two consecutive of which differ by replacing half the boundary of some polygon of  $|W|_0$ by the
other half. Since $X$ is a CAT(0) convex subcomplex it contains a polygonal face of $|W|_0$ as soon as it
contains two consecutive edges of the boundary. It follows that $X^0$ contains the vertices of any combinatorial
geodesic joining two of its points.

For any edge $e$ with endpoints $x\in X,y\not\in X$ we claim that the geometric wall $m$ separating $x$ from $y$
does not separate $x$ from any other vertex $z$ of $X$. Indeed  any vertex separated from $x$ by $m$ can be
joined to $x$ by a combinatorial geodesic through $y$. So by combinatorial convexity $X$ would contain  $y$,
contradiction. This shows that $X$ is contained in the intersection $\tilde X$ of closed half-spaces whose
boundary wall separates an edge with one endpoint in $X$ and the other one outside.

We claim that $\tilde X$ contains no vertex outside $X$. Indeed let $v\not\in X^0$ denote some vertex. Choose a
vertex $w\in X^0$ such that the combinatorial distance $d(v,w)$ is minimal. Consider any geodesic from  $w$ to
$v$. Then the first edge $e$ of this geodesic ends at a vertex $y\not\in X$, and the wall separating $w$ from
$y$ does not separate $y$ from $v$. Thus $v\not\in \tilde X$. Since $X^0\subset \bar X\subset \tilde X$ and
$\tilde X$ is the union of chambers with center in $X^0$, it follows that $\bar X= \tilde X$. Since $\Gamma$ is
cofinite on $X^0$ by assumption  it follows that $\Gamma$ is cocompact on $\bar X$, and we may apply
Theorem~\ref{thmA}.
\end{proof}

The \textbf{algebraic flat rank} of a group $\Gamma$, denoted
$\algrk(\Gamma)$, is the maximal $\Z$-rank of abelian subgroups of
$\Gamma$. The \textbf{geometric flat rank} of a ${\rm CAT}(0)$
space $X$, denoted $\rk(X)$, is the maximal dimension of
isometrically embedded flats in $X$. As an immediate consequence
of Theorem~\ref{thmA} combined with the flat torus theorem, one
obtains:

\begin{corintro}\label{corB}
Let $X$ and $\Gamma$ be as in Theorem~\ref{thmA}. Then
$\rk(X)=\algrk(\Gamma)$. In particular, one has
$\rk(|W|_0)=\algrk(W)$.
\end{corintro}

It is an important result of Daan Krammer \cite[Theorem
6.8.3]{Kr94} that the algebraic flat rank of $W$ can be easily
computed in the Coxeter diagram of $(W,S)$.

\smallskip
The equality between the algebraic flat rank of $W$ and the
geometric flat rank of $|W|_0$ was conjectured in \cite{BRW}.
Actually, it is shown in loc.~cit. that this equality allows to
compute very efficiently the so called \emph{(topological) flat
rank} of certain automorphism groups of locally finite buildings
whose Weyl group is $W$. The groups in question carry a canonical
structure of locally compact totally discontinuous topological
groups; furthermore they are topologically simple \cite{Re04}. The
topological flat rank mentioned above is a natural invariant of
the structure of topological group (see \cite{BRW} for more
details).

\smallskip
The class of pairs $(X,\Gamma)$ satisfying the assumptions of Theorem~\ref{thmA} is larger than one might
expect. Assume for example that $\Gamma$ acts geometrically by cellular isometries on a CAT(0) cubical complex
$X$, and that $\Gamma$ acts in a `special' way on hyperplanes:
\begin{enumerate}
\item for any hyperplane $H$ of $X$ and any element $g\in\Gamma$, either $gH=H$, or $H$ and $gH$ have disjoint
neighbourhoods \item for any two distinct, intersecting hyperplanes $H,H'$ of $X$ and any element $g\in\Gamma$,
either $gH'$ intersects $H$, or $H$ and $gH'$ have disjoint neighbourhoods
\end{enumerate}
Such `special' actions are studied in \cite{HW06}, where it is proved that in the above situation there exists a
right-angled Coxeter group $W$, an embedding $\Gamma\to W$ and an equivariant cellular isometric embedding $X\to |W|_0$. Thus Corollary~\ref{thmA'} applies to groups acting geometrically and specially on  CAT(0) cubical complexes.
When the action is free we obtain:

\begin{corintro}\label{corC}
The fundamental group of a compact non positively curved special cube complex is hyperbolic iff it does not contain $\Z\times \Z $.
\end{corintro}

The fundamental groups of the ``clean'' ($VH$-)square complexes studied in  \cite{Wi05}  are examples of
virtually special groups (by Theorem 5.7 of \cite{HW06}). Thus our Theorem~\ref{thmA}  provides in this case a
new proof of the equivalence between hyperbolicity and absence of $\Z\times \Z $. Note that Wise's result
applies to malnormal or cyclonormal $VH$-complexes, which are a priori more general than the virtually clean
ones. But in \cite{Wi05} Wise asks explicitly wether malnormal or cyclonormal implies virtually clean; and he
proved already this converse implication for many classes of $VH$-complexes.

\smallskip
Not surprisingly, Theorem~\ref{thmA} also provides a control on
geometric flats isometrically embedded in the ${\rm CAT}(0)$
realization of arbitrary Tits buildings. More precisely, we have:

\begin{theointro}\label{thmC}
Let $(W,S)$ be a Coxeter system and $\mathscr{B}$ be a building of type $(W,S)$. Every geometric flat of the
${\rm CAT}(0)$ realization $| {B}|_0$ of $ {B}$ is contained in an apartment. In particular, one has $\rk(|
{B}|_0)=\algrk(W)$.
\end{theointro}
Note that in  \cite{BRW} the authors had established the equality $\rk(|
{B}|_0)=\rk(|
{W}|_0)$.

Finally, we recall from \cite[Theorem~B]{Kl99} that if $X$ is a
locally compact complete ${\rm CAT}(0)$ space on which $\Isom(X)$
acts cocompactly, then the geometric flat rank of $X$ coincides
with five other quantities, among which the following ones:
\begin{itemize}
\item The maximal dimension of a quasi--flat of $X$.

\item $\sup\{k \; | \; H_{k-1}(\partial_T X) \neq \{0\}\}$, where
$\partial_T X$ denotes the Tits boundary of $X$.

\item The geometric dimension of any asymptotic cone of $X$.
\end{itemize}

This applies of course to the Davis complex $|W|_0$, but also to
many locally finite buildings of arbitrary type, including all
locally finite Kac-Moody buildings. In particular,
Corollary~\ref{corB} and Theorem~\ref{thmC} above, combined with
Daan Krammer's computation of $\algrk(W)$, provide a very
efficient way to compute all these quantities for these examples.


In Section~\ref{section:preliminaries}, we first recall basic facts on the Davis--Moussong geometric realization
of Coxeter groups. In particular we introduce the walls, the half-spaces and the chambers.

In Section~\ref{section:combinatorialconvexity} we define  combinatorial convex subsets of the Davis--Moussong
geometric realization, and we establish an important Lemma.

In Section~\ref{section:eucltriangle} we present the main technical tools of this article. If a family of walls
behaves as if it was contained in a Euclidean triangle subgroup, then in fact it generates a Euclidean triangle
subgroup (see Lemmas~\ref{lem:triangle1} and~\ref{cor:triangle2}  for precise statements).

In Section~\ref{section:wallflat} we describe completely the combinatorial structure of the set of walls
separating a given flat. The reflections along these walls generate a subgroup that we also describe.

In Section~\ref{section:coxrank} we explain how to get a rank $n$ free abelian group out of a rank $n$ flat.

And  in Section~\ref{section:buildrank} we explain how to deduce the statement on buildings from the statement
on Coxeter complexes.

\section{Preliminaries}\label{section:preliminaries}
Let $(W,S)$ be a Coxeter system with $S$ finite. The Davis complex
associated with $(W,S)$, denoted $|W|_0$, is a $\mathrm{CAT}(0)$
cellular complex equipped with a faithful, properly discontinuous,
cocompact action of $W$ (see \cite{Da98}).

Recall that a \textbf{reflection} of $W$ is, by definition, any
conjugate of an element of $S$. The fixed point set of a
reflection in $|W|_0$ is called a \textbf{wall}. Note that a wall
is a closed convex subset of $|W|_0$. A fundamental property is
that every wall separates $|W|_0$ into two open convex subsets,
whose respective closures are called \textbf{half-spaces}. If $a$
is a half-space, its boundary is a wall which is denoted by
$\partial a$. If $x \in |W|_0$ is a point which is not contained
in any wall, then the intersection of all half-spaces containing
$x$ is compact; this compact set is called a \textbf{chamber} of
$|W|_0$. The $W$-action on the chambers of $|W|_0$ is free and
transitive.

Let $x,y$ denote two non-empty convex subsets of $|W|_0$. We say that a wall $m$ separates $x$ from $y$ whenever
$x$ is contained in one of the half-spaces delimited by $m$, $y$ is contained in the other half-space, and
neither $x$ nor $y$ are contained in $m$.


We will use the following notation. Given a wall $m$ of $|W|_0$,
the unique reflection fixing $m$ pointwise is denoted by $r_m$. For any
set $M$ of walls, we set $W(M):=\la r_m | \; m \in M \ra$.
Recall that $W(M)$ is itself a Coxeter system on a certain set of reflections $(r_\nu)_{\nu\in N}$, where each
wall $\nu\in N$ is of the form $\nu=w\mu$ for some $w\in W(M)$ and some $\mu\in M$ (see \cite{Deo89}). Such a
subgroup will be called a \textbf{reflection subgroup}.

Finally, given two points (resp. two convex subsets) $x, y$ of $|W|_0$, we denote by $ \mathscr{M}(x,y)$ the set
of all walls which separate $x$ from $y$.  Two chambers $c,c$ are said to be \textbf{adjacent} whenever
$\mathscr{M}(c,c')$ is empty, or consists in a single wall $m$ (in which case $r_m(c)=c'$). A gallery (of length
$n$) is a sequence $(c_0,c_1,\dots,c_n)$ of chambers such that $c_i$ and $c_{i+1}$ are adjacent chambers for
$i=0,\dots,n-1$. The gallery defines a unique sequence of walls it crosses (this sequence might be empty if the
gallery is a constant sequence).

We get a (discrete) distance on the set of chambers by considering the infimum of the length of all galleries
from the first chamber to the second. Using the simple transitive action of $W$ on the chambers, this gallery
distance is identified with the word metric on $(W,S)$.

It is well known that for two chambers $c,c'$ the gallery distance $d_{\rm gal}(c,c')$ is the cardinality of
$\mathscr{M}(c,c')$, and that a gallery from $c$ to $c'$ has length $d_{\rm gal}(c,c')$ if and only if the
sequence of walls it crosses has no repetition. Furthermore for any gallery from $c$ to $c'$ the set of walls
separating $c$ from $c'$ is the set of walls appearing an odd number of times in the sequence of walls that the
gallery crosses.

The following basic lemmas are  well known; their proofs are easy exercises.

\begin{lem}\label{lem:prebasic-Coxeter}
Let  $x,y$ be two points of $|W|_0$. There are two  chambers $c_x$, $c_y$  such that $x \in c_x$, $y \in c_y$
and $ \mathscr{M}(x,y)= \mathscr{M}(c_x, c_y)$.\qed
\end{lem}

\begin{lem}\label{lem:basic-Coxeter}
Let $x, y \in |W|_0$. There exists $\gamma \in W( \mathscr{M}(x, y))$ such that $x$ and $\gamma.y$ are contained
in a common chamber.\qed
\end{lem}

\section{Combinatorial convexity}\label{section:combinatorialconvexity}

A subset $F \subset |W|_0$ is called \textbf{combinatorially convex} if either $F = |W|_0$ or $F$ coincides with
the intersection of all half-spaces containing it. The \textbf{combinatorial convex closure} of a subset $F
\subset |W|_0$ will be denoted by ${\rm Conv}(F)$. Hence ${\rm Conv}(F)$ is either the whole $|W|_0$ (if $F$ is
not contained in any half-space) or the intersection of all half-spaces of $|W|_0$ containing $F$. Since
half-spaces are subcomplexes of the first barycentric subdivision of $|W|_0$ we note that combinatorially convex
subsets are subcomplexes as well.

Since half-spaces are CAT(0) convex, combinatorially convex subcomplexes are CAT(0) convex, but we will rather
use the following elementary combinatorial convexity property:
 all chambers of a geodesic gallery from a chamber $c$ to a chamber $c'$ belong to ${\rm Conv}(c\cup c')$.

\begin{lem}\label{lem:convexity}
Let $x, y \in |W|_0$ and assume that the set $ \mathscr{M}(x,y)$ possesses a subset $M$ such that for all $m \in
M$ and $\mu \in \overline M= \mathscr{M}(x,y) \backslash M$, the reflections $r_m$ and $r_\mu$ commute. Then the
combinatorial convex closure of $\{x, y\}$ contains a point $z$ such that $ \mathscr{M}(y, z)=M$ and
$\mathscr{M}(x,z)=\overline M$.
\end{lem}
\begin{proof}
Let $c_x$, $c_y$ be chambers such that $x \in c_x$, $y \in c_y$ and $ \mathscr{M}(x,y)= \mathscr{M}(c_x, c_y)$
(see Lemma~\ref{lem:prebasic-Coxeter}). We prove that there exists a chamber $c_z$ such that $ \mathscr{M}(c_y,
c_z)=M$ and $ \mathscr{M}(c_x, c_z)=\overline M$ (note that such a chamber necessarily lies  in the
combinatorial convex closure of $c_x\cup  c_y$).

This implies the desired result. Indeed since $\mathscr{M}(x,y)= \mathscr{M}(c_x, c_y)$ we have ${\rm
Conv}(\{x,y\})={\rm Conv}(c_x\cup c_y)$. Furthermore since $\mathscr{M}(y,c_y)=\emptyset$ we have
$\mathscr{M}(c_z,y)\subset \mathscr{M}(c_z,c_y)$. Conversely if $m\in \mathscr{M}(c_z,c_y)$ then $m$ does not
separate $c_z$ from $c_x$ -- otherwise $c_z$ would not be inside ${\rm Conv}(c_x\cup c_y)$. Thus $m$ separates
$c_y$ from $c_x$, and so $m\in \mathscr{M}(x,y)$. In particular $y\not\in m$. Thus in fact $m\in
\mathscr{M}(c_z,y)$. Consequently $\mathscr{M}(c_z,y)= \mathscr{M}(c_z,c_y)\ (=M)$, and similarly
$\mathscr{M}(c_z,x)= \mathscr{M}(c_z,c_x)\ (=\overline M)$. We then define $z$ to be any point in the interior
of the chamber $c_z$.

It remains to prove the statement for chambers. To this end, we argue  by induction on the cardinality $n$ of $
\mathscr{M}(c_x,c_y)$. We may assume $n>0$.

Consider some geodesic gallery $(c_0=c_x,\dots,c_{n-1},c_n=c_y)$. Let $\mu$ denote the unique wall separating
$c_{n-1}$ from $c_n$.  By induction there is a chamber  $d$  such that $\mathscr{M}(c_x,d)=\overline
M\setminus\{\mu\}, \mathscr{M}(d,c_{n-1})=M\setminus\{\mu\}$. We then have $\mathscr{M}(d,c_y) =
\mathscr{M}(d,c_{n-1})\cup\{\mu\}$.

If $\mu\in M$, then  the chamber $d$ satisfies $\mathscr{M}(c_x,d)=\overline M$ and $\mathscr{M}(d,c_y) = M$, so
we are done.

Assume now that $\mu\in\overline M$, so $M=\mathscr{M}(d,c_{n-1})$. Consider a gallery  from $d$ to $c_{n-1}$ of
minimal length. If this gallery has length 0 then $M=\emptyset$ and we take $c_z=c_y$. Otherwise let $m\in M$
denote the last wall that the gallery crosses. Let $d'$ denote the chamber $r_mr_\mu(c_{n-1})$. Then $d'$ is
adjacent to $c_{n-2}$, and $d'$ is also adjacent to $c_n$ because $r_mr_\mu=r_\mu r_m$. It follows that there
exists a gallery of minimal length from $c_x$ to $c_y$ whose last crossed wall is $m$. So in fact we are back to
the first case, and thus we are done.
\end{proof}

Note that the corresponding statement (for vertices) is true in an arbitrary CAT(0) cubical complex $X$. Indeed
for any two vertices $x,y$ of $X$ such that the set $ \mathscr{M}(x,y)$ of hyperplanes of $X$ separating $x$
from $y$ may be written $ \mathscr{M}(x,y)=M\sqcup \overline M$ so that every hyperplane of $M$ is perpendicular
to every hyperplane of $\overline M$, there exists a vertex $z$ such that $ \mathscr{M}(z,y)=M$ and
$\mathscr{M}(z,x)=\overline M$. Clearly $z$ is on some combinatorial geodesic from $x$ to $y$, thus $z$ is in
the convex hull of $\{x,y\}$.

\section{The Euclidean triangle lemmas}\label{section:eucltriangle}

In what follows, a \textbf{Euclidean triangle subgroup} of the Coxeter group $W$ is a reflection subgroup which
is isomorphic to one of the three possible irreducible Coxeter groups containing $\Z\times \Z$ as a finite index
subgroup. We say that a set $P$ of walls is  \textbf{Euclidean} whenever there exists a wall $m$ such that
$P\cup\{m\}$ generates a Euclidean triangle subgroup of $W$. We will be mainly interested in the case when $P$
is a set of pairwise disjoint walls.

The following lemma relates the combinatorial configuration of a certain set of walls $M$ of $|W|_0$ with the
algebraic structure of $W(M)$. This provides the key ingredient which allows to understand the walls of a
geometric flat of $|W|_0$, see Proposition~\ref{prop:flat-walls} below.



\begin{lem}\label{lem:triangle1}
There exists a constant $L$, depending only on the Coxeter system $(W,S)$, such that the following property
holds. Let $a, b, h_0, h_1, \dots, h_n$ be a collection of half-spaces of $|W|_0$ such that:
\begin{itemize}
\item[(1)] $\emptyset \neq a \cap b \subsetneq h_0 \subsetneq h_1 \subsetneq \dots \subsetneq h_n$,

\item[(2)] $\emptyset \neq \partial a \cap \partial b \subset \partial h_0$,

\item[(3)] $\partial a$ and $\partial b$ both meet $\partial h_i$ for each $i = 1, \dots, n$.
\end{itemize}
If $n \geq L$, then the group generated by the reflections through the walls $\partial a$, $\partial b$,
$\partial h_0$, $\partial h_1, \dots, \partial h_n$ is  a {Euclidean triangle subgroup}.
\end{lem}
\begin{proof}
See \cite[Theorem~A]{Cap}.
\end{proof}

A set $P$ of walls of $|W|_0$ is called a \textbf{chain of walls} if there exists a set $A$ of half-spaces  of
$|W|_0$ such that $A$ is totally ordered by inclusion and $P=\{\partial a,a\in A\}$ (for short we write
$P=\partial A$). There are three kinds of chains of walls. We say that $P$ is a \textbf{segment of walls} if it
is a finite chain of walls. We say that  $P$ is a \textbf{line of walls} if $P=\partial A$, with  $A$ a set of
half-spaces such that the ordered set $(A,\subset)$ is isomorphic to $(\integers,\le)$. And we say that  $P$ is
a \textbf{ray of walls} if $P=\partial A$, with  $A$ a set of half-spaces such that the ordered set
$(A,\subset)$ is isomorphic to $(\naturals,\le)$.

\begin{lem}\label{lem:propeucl} Let $P$ denote a nonempty set of walls which are all disjoint from a given
wall $\mu$.  Assume that $P\cup\{\mu\}$ is Euclidean. Then $P\cup\{\mu\}$ is a chain and $W(P\cup\{\mu\})$ is
infinite dihedral.
\end{lem}
\begin{proof} Let $\mu'$ denote some wall such that $W(P \cup \{\mu,\mu'\})$ is a Euclidean triangle subgroup.
Represent $W(P \cup \{\mu,\mu'\})$ as a group of isometries of the Euclidean plane (in such a way that the
abstract reflections  act as geometric reflections).

Let $m,m'$ denote two walls of $P\cup\{\mu\}$. Note that $m\cap m'=\emptyset$ if and only if the order of
$r_mr_{m'}$ is infinite. In the geometric representation we have $m\cap m'=\emptyset$ if and only if the
Euclidean lines $L(m), L(m')$ fixed pointwise by $m$  and $m'$ are parallel. Since we assume $m\cap
\mu=\emptyset$ or $m=\mu$, we deduce that $L(m)$ is parallel  to $L(\mu)$. Similarly $L(m')$ is parallel to
$L(\mu)$. Thus $L(m)$ and $L(m')$ are parallel, which implies that $m=m'$ or $m\cap m'=\emptyset$.

Thus $P \cup \{\mu\}$ is a set of pairwise disjoint walls (of cardinality $\ge 2$). By looking at the geometric
representation we deduce that $W(P\cup\{\mu\})$ is infinite dihedral.
Note  that the set of walls associated with all the reflections of any infinite dihedral reflection subgroup is
a line of walls (this can be seen by considering a generating set consisting of two reflections; the associated
walls cut $|W|_0$ into three pieces, one of which is a fundamental domain for the reflection subgroup that we
consider). It follows that $P \cup \{\mu\}$ is a chain.
\end{proof}

Let $T$ denote any subset of the generating set $S$. Then any conjugate of the subgroup $W(T)$ is called  a
\textbf{parabolic subgroup}. The \textbf{parabolic closure} of any subgroup $\Gamma\subset W$ is the
intersection of all parabolic subgroups of $W$ containing $\Gamma$; we denote it by $\widetilde\Gamma$. With
this terminology, we have:
\begin{lem}\label{lem:triangle2}
Let $P$ be a set of pairwise disjoint walls of $|W|_0$. Assume that there exists a wall $m$ such that $W(P \cup \{m\})$ is a
Euclidean triangle subgroup. Then the parabolic closure $\widetilde{W(P)}$  satisfies the
following conditions:
\begin{enumerate}
\item $\widetilde{W(P)}$ is isomorphic to an irreducible affine Coxeter group.

\item For all walls $\mu, \mu', \mu''$, if $\mu$ separates $\mu'$ from $\mu''$ and if $r_{\mu'}$ and $r_{\mu''}$
both belong to $\widetilde{W(P)}$, then $r_\mu$ also belongs to $\widetilde{W(P)}$.

\item For any line of walls $P'$ and any wall $\mu$, if $W(P') \leq \widetilde{W(P)}$ and if $W(P' \cup
\{\mu\})$ is a Euclidean triangle subgroup, then $r_\mu$ belongs to $\widetilde{W(P)}$.
\end{enumerate}
\end{lem}
\begin{proof}
Point (1) follows from a theorem of D.~Krammer which appears in \cite[Theorem~1.2]{CM05} (see also Theorem~3.3
in loc. cit.); (2) and (3) follow from (1) using convexity arguments, see \cite[Lemma~8]{Cap} for details.
\end{proof}
We may now deduce an other useful result of the same kind as Lemma~\ref{lem:triangle1}:

\begin{cor}\label{cor:triangle2}
Let $P$ be a set of pairwise disjoint walls of $|W|_0$ and let $m$ be a wall such that $W(P \cup \{m\})$ is a
Euclidean triangle subgroup. Then $W$ possesses a Euclidean triangle subgroup, denoted by $\overline{W(P \cup
\{m\})}$, containing $W(P \cup \{m\})$ and such that $r_\mu \in \overline{W(P \cup \{m\})}$ for each wall $\mu$
satisfying either of the following conditions:
\begin{enumerate}
\item There exist $\mu', \mu'' \in P$ such that $\mu$ separates $\mu'$ from $\mu''$.

\item $\mu$ is disjoint from $m$ and moreover $W(P \cup \{\mu\})$ is a Euclidean triangle subgroup.
\end{enumerate}
\end{cor}
\begin{proof}
Let $\widetilde{W(P)} \leq W$ be the irreducible affine Coxeter group provided by Lemma~\ref{lem:triangle2}. By
Lemma~\ref{lem:triangle2}(3) we have $r_m \in \widetilde{W(P)}$. Let $P'$ be the set consisting of all those
walls $p'$ such that $r_{p'} \in \widetilde{W(P)}$ and that there exists $p \in P \cup \{m\}$ which does not
meet $p'$. Define $\overline{W(P \cup \{m\})}:=W(P' \cup P \cup \{m\})$. The group $\overline{W(P \cup \{m\})}$
is a Euclidean triangle subgroup, because it is a subgroup of an affine Coxeter group generated by reflections
corresponding to two directions of hyperplanes. Given a wall $\mu$ satisfying (1) or (2), we obtain successively
$r_\mu \in \widetilde{W(P)} $ by Lemma~\ref{lem:triangle2} and then $\mu \in P'$ by the definition of $P'$.
\end{proof}

\section{The walls of a geometric flat}\label{section:wallflat}

Let $F$ be a geometric flat which is isometrically embedded in the Davis complex $|W|_0$ of $W$. Let $
\mathscr{M}(F)$ denote the set of all walls which separate points of $F$:
$$ \mathscr{M}(F):= \bigcup_{x, y \in F}  \mathscr{M}(x,y).$$

\begin{lem}\label{0}
For all $\mu \in  \mathscr{M}(F)$, the set $\mu \cap F$ is a Euclidean hyperplane of $F$.
\end{lem}
\begin{proof}
Let $x , y$ be points of $F$ which are separated by $\mu$. We know
that $\mu \cap F$ is a closed convex subset of $F$ which separates
$F$ into two open convex subsets. Thus the result will follow if
we prove that the geodesic segment $[x, y]$ joining $x$ to $y$
meets $\mu$ in a single point. This is a local property, which can
easily be  checked in a single (Euclidean) cell of $|W|_0$ (see
\cite[Lemma~3.4]{NV02} for details).
\end{proof}

\begin{lem}\label{lem:wall-neighborhood}
Let $\mu$ be a wall which meets $F$. Assume that $F$ contains a  Euclidean half-space  $F^+$ such that $F^+
\cap \mu \neq\emptyset$ and $F^+$ is contained in a $\varepsilon$--neighborhood of $\mu$ for some $\varepsilon
> 0$. Then $F \subset \mu$.
\end{lem}
\begin{proof}
Let $d$ be the distance function of the Davis complex $|W|_0$. Since $\mu$ is a closed convex subset, the
function $d_\mu: |W|_0 \to \mathbf{R}^+: x \mapsto \inf\{d(x,y) | \; y \in \mu\}$ is convex (see \cite[\S
II.2]{BH99}). By assumption, the restriction $d_\mu|_{F^+}$ of $d_\mu$ to $F^+$ is bounded. Therefore
$d_\mu|_{F^+}$ must be constant, as it is the case for any bounded convex function on an unbounded convex
domain. Since $\mu$ meets $F^+$ by hypothesis, we have $d_\mu|_{F^+} = 0$, that is to say, $F^+ \subset \mu$. By
Lemma~\ref{0}, this implies $F \subset \mu$.
\end{proof}

Two elements $\mu, \mu' $ of $ \mathscr{M}(F)$ will be called \textbf{$F$--parallel} if their respective traces
on $F$ are parallel in the Euclidean sense. In symbols, this writes: $$\mu \|_F \;  \mu' \hspace{.5cm}
\Leftrightarrow \hspace{.5cm} \mu \cap F = \mu' \cap F \text{ or } \mu \cap F \cap \mu' = \emptyset.$$ The
relation of $F$--parallelism is an equivalence relation on $ \mathscr{M}(F)$.

Besides the relation of $F$--parallelism, there is an other relation of \textbf{global parallelism} on the walls
of $F$ defined by
$$\mu \| \mu' \hspace{.5cm} \Leftrightarrow \hspace{.5cm} \mu =
\mu' \text{ or } \mu \cap \mu' = \emptyset.$$ Clearly $\mu \| \mu'  \Rightarrow \mu \|_F \;  \mu'$. Given $\mu
\in  \mathscr{M}(F)$, we set $P_F(\mu):=\{ m \in  \mathscr{M}(F) \; | \; m \| \mu \}$. Thus $P_F(\mu)$ is
contained in the $F$--parallel  class of $\mu$. Note that, in contrast with the $F$--parallelism, the relation
of global parallelism is not transitive in general:  two distinct walls of $P_F(\mu)$ may have non trivial
intersection.

Any large set of walls contains two non-intersecting ones (see \cite[Lemma~3]{NR03}). Consequently, the set of
$F$--parallel  classes is finite. Since chambers are compact and $F$ is unbounded, it follows that some
$F$--parallel  class must be infinite. Actually, all of them are, as follows from the following:

\begin{lem}\label{lem:infinite--parallel-class}
Given any $\mu \in  \mathscr{M}(F)$, there exist two rays of walls $M^+(\mu), M^-(\mu) \subset \mathscr{M}(F)$
such that $\mu$ separates any element of $M^+(\mu)$ from any element of $M^-(\mu)$. In particular, $\mu$ does
not meet any element of $M^+(\mu) \cup M^-(\mu)$, and $P_F(\mu)$ contains a line of walls (passing through
$\mu$).
\end{lem}
\begin{proof} Consider a line of $F$ which meets orthogonally the $F$-hyperplane $\mu\cap F$. Using
Lemma~\ref{lem:wall-neighborhood} we see that when a point $p$ goes at infinity on the line, its distance to
$\mu$ must tend to infinity. Now by the so called \emph{parallel wall theorem} (see \cite[Theorem~2.8]{BH94})
any point at large distance from a given wall in $|W|_0$ is separated from that wall by some other wall of
$|W|_0$. The Lemma follows.
\end{proof}

\begin{rem}\label{rem:chain}  For $\mu\in\MF$, any subset $P\subset P_F(\mu)$ of
pairwise disjoint walls is a chain of walls. Indeed for three distinct walls $p_1,p_2,p_3\in P$ we have
$p_i\|_F\mu$, thus  $p_1,p_2,p_3$ are mutually $F$--parallel. The Euclidean hyperplanes $p_i\cap F$ are pairwise
disjoint, so we may assume that  $p_2\cap F$ separates $p_1\cap F$ from $p_3\cap F$. It follows that $p_2$
separates $p_1$ from $p_3$. Hence $\{p_1,p_2,p_3\}$ is a segment of walls. Since any $3$--subset of $P$ is a
chain, it follows that $P$ itself is a chain.
\end{rem}

We will see in Proposition~\ref{prop:groupe-Meucl} below
that the restriction of the relation of global parallelism to a
certain subset $\Meucl$ of $\MF$ is an equivalence.

By definition, the subset $\Meucl \subset \MF$ consists of all
those walls $\mu \in \MF$ which satisfy the following property:

\smallskip
\noindent \emph{There exists a wall $\mu' \in \MF$ such that
$W(P_F(\mu) \cup \{\mu'\})$ is a Euclidean triangle subgroup. }

\smallskip Applying Lemmas~\ref{lem:propeucl} and~\ref{lem:infinite--parallel-class}, we get the following:

\begin{lem}\label{lem:propertymeucl}
Assume $\mu\in\Meucl$, and more precisely that $W(P_F(\mu)\cup\{\mu'\})$ is a Euclidean triangle subgroup for
some $\mu'\in\MF$. Then:
\begin{enumerate}

\item[(i)] $P_F(\mu)$ is a line of walls.

\item[(ii)] For all $m \in P_F(\mu)$, one has $P_F(\mu) \subset P_F(m)$. In particular  $P_F(\mu)=P_F(m)$
provided $m\in\Meucl$.

\item[(iii)] $W(P_F(\mu))$ is an infinite dihedral subgroup of $W$, and is a maximal one.

\item[(iv)] $r_{\mu'}$ does not centralize $W(P_F(\mu))$.\qed
\end{enumerate}
\end{lem}

The following
lemma outlines the main combinatorial properties of the set
$\Meucl$.

\begin{lem}\label{lem:Meucl} We have the following:
\begin{itemize}
\item[(i)] Let $P \subset \MF$ be a line of walls. If there exists $m \in \MF$ such that the group $W(P \cup
\{m\})$ is a  Euclidean triangle subgroup, then $P \subset \Meucl$.

\item[(ii)] Let $m \in \MF$. If $m \not \in \Meucl$, then $m$
meets every element of $\Meucl$.

\item[(iii)] Let $m, m' \in \MF$. If the reflections $r_m$ and $r_{m'}$ do not commute and if $m$ and $m'$ are
not $F$--parallel, then $m \in \Meucl$.

\item[(iv)] Let $m, m' \in \MF$. If the reflections $r_m$ and
$r_{m'}$ do not commute and if $m' \in \Meucl$, then $m \in
\Meucl$.
\end{itemize}
\end{lem}

Before proving the lemma, it is convenient to introduce the following additional terminology. A set $P$ of walls
of $|W|_0$   is said to be \textbf{convex} whenever the following holds:  for each wall $m$ of $|W|_0$
separating two walls of $P$, we have $m \in P$. For example, for all $x, y \in |W|_0$ the set $ \mathscr{M}(x, y
)$ is convex; moreover, the set $\MF$ is convex as well.

\begin{proof}[Proof of Lemma~\ref{lem:Meucl}(i)]
Let $\mu \in P$. Since $P \subset \MF$  is a line of walls we have $P \subset P_F(\mu)$. There are finitely many
walls separating two disjoint walls of $|W|_0$. The line of walls $P$ may be written as a union of segments of
walls $\{\mu_n,\mu_{n+1}\}$ ($n\in\integers$) so that no $m\in P$ separates $\mu_n$ from $\mu_{n+1}$. Choose
then a segment of walls $P_n \subset P_F(\mu)$ such that  $P_n \cap P =\{\mu_n,\mu_{n+1}\}$ and any wall $m\in
P_n\setminus \{\mu_n,\mu_{n+1}\}$ separates $\mu_n$ from $\mu_{n+1}$ and moreover $P_n$ is maximal with respect
to these properties. Set $\bar{P}=\cup_k P_k$. Then $P\subset \bar{P}\subset P_F(\mu)$, $\bar{P}$ is a line and
for every wall $m'$ of $P_F(\mu)\setminus \bar P$ the set $\bar P\cup\{m'\}$ is not a line anymore.

By construction for every $p \in \bar{P}$ there exist $p', p'' \in P$ such that $p$ separates $p'$ from $p''$.
Therefore, since $W(P \cup \{m\})$ is a Euclidean triangle subgroup, we have $W(\bar{P} \cup \{m\}) \subset
\overline{W(P \cup \{m\})}$ by Corollary~\ref{cor:triangle2}. In particular, $W(\bar{P} \cup \{m\})$ is a
Euclidean triangle subgroup. Hence we are done if we show that $\bar{P} = P_F(\mu)$. This is what we do now.

\smallskip Let $m'$ denote a wall separating two walls $p',p''$ of $\bar P$. Then $m'\in\MF$ and
by Corollary~\ref{cor:triangle2} the subset $\bar P\cup\{m'\}$ is still Euclidean. By Lemma~\ref{lem:propeucl}
$\bar P\cup\{m'\}$ is a line, and by the maximality of $\bar{P}$ we have $m'\in \bar P$. Thus $\bar{P}$ is a
convex set of walls.

Assume by contradiction that there exists $m' \in
P_F(\mu) \backslash \bar{P}$. By the maximality of $\bar{P}$, the
set $\bar{P}\cup\{m'\}$ is not a line anymore.  By Remark~\ref{rem:chain} this implies that $m'$
 meets at least one element of $\bar{P}$.
Let $\bar{P}'$ denotes the (nonempty) subset of $\bar{P}$ consisting of all those walls which
meet $m'$.  Note that by the
definition of $\bar{P}'$, for all $p \in \bar{P}$, if there exist
$p', p'' \in \bar{P}'$ such that $p$ separates $p'$ from $p''$,
then $p \in \bar{P}'$. Since $\bar{P}$ is convex, this shows in
particular that $\bar{P}'$ is convex.

If $\bar{P}'$ is finite, it is a segment of the line $\bar P$ and  there exist $p', p'' \in \bar{P}$
such that $m'$ separates $p'$ from $p''$. Since $\bar{P}$ is
convex,  this implies that $m' \in \bar{P}$, a contradiction.

Hence $\bar{P}'$ is infinite. Since $\mu\not\in \bar{P}'$ and $ \bar{P}'$ is convex, we see that $ \bar{P}'$ is
a ray of walls (contained in $\bar P$, and not containing $\mu$).

By  Lemma~\ref{lem:propeucl} the group $W(\bar{P})$ is infinite dihedral. Since $\bar{P}$ is a line of walls,
the wall $\pi$ of any reflection $r_\pi$ of  $W(\bar{P})$ separates two walls $p',p''$ of $\bar{P}$. By
convexity we then have $\pi\in\bar{P}$: the reflections of $W(\bar{P})$ are precisely the reflections along
walls of $\bar P$. We note two consequences of that. Firstly $\bar P$ is invariant under $W(\bar P)$.  Secondly
we have $W(\bar{P})=W(\bar{P}_0)$ for any convex subset $\bar{P}_0 \subset \bar{P}$ of cardinality at least $2$.
In particular we have $W(\bar{P})=W(\bar{P}')$.

The reflection $r_{m'}$ does not centralize $W(\bar{P}')$, otherwise it would centralize $W(\bar{P})$ and,
hence, $m'$ would meet $\mu$. Consequently  $r_{m'}$ does not centralize $W(\bar{P}'_0)$ for all convex subset
$\bar{P}'_0 \subset \bar{P}'$ of cardinality at least~$2$. Hence there are infinitely many walls $\bar{p}'$ in
the ray $\bar{P}'$ such that the reflections $r_{m'}$ and $r_p$ do not commute. Let $\bar{p}'\in \bar{P}'$
denote some wall such that the reflections $r_{m'}$ and $r_{\bar{p}'}$ do not commute, and that the collection
of all walls of $\bar{P}'$ which separate $\bar{p}'$ from $\mu$ is of  cardinality greater than  the constant
$L(\ge 1)$ of Lemma~\ref{lem:triangle1}.

Let $m'':=r_{\bar{p}'}(m')$. Let $\{\bar{p}_1,\dots,\bar{p}_k\}$ denote the segment of walls of $\bar P'$ which
separate $\mu$ from $\bar{p}'$ (we have $k\ge L$). Then the walls $\bar{m}_i=r_{\bar{p}'}(\bar{p}_i)$ belong to
the ray $\bar{P}'$ by convexity (remember that $r_{\bar{p}'}(\mu)\in \bar P$). Hence each of them meets $m'$. By
construction each of them also meets $m''$. By Lemma~\ref{lem:triangle1} we deduce that
$W(m',m'',\bar{m}_1,\dots,\bar{m}_k,\bar{p}') =W(m',\bar{m}_1,\dots,\bar{m}_k,\bar{p}') $ is a Euclidean
triangle subgroup. Since $\{\bar{m}_1,\dots,\bar{m}_k,\bar{p}'\}$ is a convex subsegment of $\bar P$ containing
at least two walls we see that $W(\bar{P} \cup \{m'\})$ is a Euclidean triangle subgroup. Since $\mu\in\bar P$
and $m'\cap\mu=\emptyset$, this contradicts Lemma~\ref{lem:propeucl}, thereby completing the proof of the
desired assertion.
\end{proof}

\begin{proof}[Proof of Lemma~\ref{lem:Meucl}(ii)]
Let $m \in \MF$. Assume that there exists $\mu \in \Meucl$ which does not meet $m$. In other words $m\in P_F(\mu)$.  By
Lemma~\ref{lem:Meucl}(i), $\mu \in \Meucl$ implies $P_F(\mu) \subset
\Meucl$. Thus $m  \in \Meucl$.

\end{proof}

\begin{proof}[Proof of Lemma~\ref{lem:Meucl}(iii)]
Let $M$ be the $F$--parallel  class of $m$ and let $m'':= r_m(m')$. Since $m$ and $m'$ are not $F$--parallel,
there are points $x,y$ on $m\cap F$ which are separated by $m'$. Thus $m''$ separates $x$ from  $y$ as well. It
follows that $m''\in \mathscr{M}(F)$.

We now show that $m''$ is not $F$--parallel to $m$. To this end, first note that $m''$ contains $m \cap m' \cap
F$ which is nonempty. Hence, if $m''$ were $F$-parallel to $m$, then we would have $m \cap F = m'' \cap F$. This
yields successively $m \cap F = r_m(m') \cap F$ and then $m \cap r_m(F) = m' \cap r_m(F)$. Since $m \cap F$ is
pointwise fixed by $r_m$, we have $m \cap F \subset m \cap r_m(F)$, whence finally $m \cap F \subset m'$, which
contradicts the fact that  $m$ and $m'$ are not $F$-parallel. This shows that $m''$ is not $F$--parallel to $m$
and it follows that $m'$ and $m''$ both meet every element of the $F$--parallel  class $M$.

By Lemma~\ref{lem:infinite--parallel-class}, $M$ contains a line $P$ containing $m$. In particular $P$ is
infinite. By Lemma~\ref{lem:triangle1}, the group $W(P \cup \{m'\})$ is a Euclidean triangle subgroup.
Therefore, we deduce from Lemma~\ref{lem:Meucl}(i) that $m \in \Meucl$.
\end{proof}

\begin{proof}[Proof of Lemma~\ref{lem:Meucl}(iv)]
Let $m \in \MF$ and ${m'} \in \Meucl$ be such that the reflections $r_m$ and $r_{m'}$ do not commute. By
Lemma~\ref{lem:Meucl}(i), we have $P_F({m'}) \subset \Meucl$. Let $P':= P_F(m')$. Hence $P'$ is a line of walls
and for all $\mu' \in P'$, we have $P_F(\mu')=P'$.

By Lemma~\ref{lem:Meucl}(ii) we may assume that $m$ meets every element of $P'$, and in fact that every element
of $P_F(m)$ meets every element of $P'$, otherwise $m \in \Meucl$ and we are done. By
Lemma~\ref{lem:infinite--parallel-class}, $P_F(m)$ contains a line of walls $P$ which contains $m$.

Let $C$ (resp. $C'$) denote the set of walls of $P$ (resp. $P'$) which meet $m'':=r_m(m')$.

Assume that $C'$ is finite. Then there exists a (convex) segment of walls $(p_+, p_1, \dots, p_n, p_-)$
contained in $P'$ such that $C'= \{p_1, \dots, p_n\}$ and $m''$ is disjoint from $p_+$ and $p_-$. We let $x_+,
x_-$ denote points lying on $m \cap p_+, m \cap p_-$ respectively. Since $m'$ separates $p_+$ from $p_-$ and $m
\setminus m'' =m\setminus m'$ we deduce that $m''$ separates $x_+$ from $x_-$. Thus $m''$ separates  $p_+$
from $p_-$. It follows that $m''\in \mathscr{M}(F)$, and in fact $m''\in P_F(p^+)$.  By hypothesis $m' \in \Meucl$, whence $P_F(p^+)=P_F(m')$. Since $m''$ meets $m'$, this implies $m' = m''$ from which it follows
that the reflections $r_m$ and $r_{m'}$ commute, a contradiction. Thus $C'$ is infinite.

By Lemma~\ref{lem:triangle1}, it follows that $W(C' \cup \{m\})$ is a Euclidean triangle subgroup. Since $P'$ is
Euclidean, we have $W(P') = W(P'_0)$ for any convex chain $P'_0 \subset P'$ of cardinality at least~$2$ (see
Lemma~\ref{lem:propeucl}). Since $C'$ is infinite and convex, we deduce  $W(P') \subset W(C' \cup \{m\})$. Since
$r_{m''}$ belongs to the Euclidean triangle subgroup  $ W(C' \cup \{m\})$ and $r_{m''}r_{\mu'}$ has finite order
for every $\mu'\in C'$, we see that  $r_{m''}r_{\mu'}$ has finite order for every $\mu'\in P'$. Thus $C'=P'$.
Moreover for all $\mu' \in P'$, the reflections $r_{\mu'}$ does not commute with $r_m$.


\smallskip  Let $\mu$ be any element of $P$ different
from $m$. Let $a$  denote the half-space bounded by $m$ such that $\mu\cap a=\emptyset$. Let $h_0$ denote the
half-space bounded by $m'$ such that $a\cap h_0\subset r_m(h_0)$. Extend $h_0$ to  a chain of half-spaces
$(h_i)_{i \in \Z}$ such that  $h_i \subset h_{i+1}$ for all $i \in \Z$ and that  $\{\partial h_i | \; i \in \Z
\} = P'$. Since $W(P'\cup\{ m\})$ is a Euclidean triangle subgroup it follows that the relation $a\cap
h_i\subset r_m(h_i)$ holds for every $i\in\Z$.  For each $i \in \Z$, choose a point $y_i \in \mu \cap \partial
h_i$  and a point $y'_i \in  \partial h_i $ in the interior of $a$. Then $y_i \in  \partial h_i $ and $y'_i \in
\partial h_i $  are separated by $m$. Since $r_m$ and $r_{\partial h_i}$ do not commute it follows that $y_i$
and $y'_i$  are separated by $r_m(\partial h_i)$. Since $y'_i\in a\cap h_i$ we deduce that  $y_i \not \in
r_m(h_i)$ for all $i \in \Z$.
 Now choose a point $x_i\in m\cap \partial  h_i$ for each $i\in\Z$.
 We have
$$x_0 \in m \cap \partial h_0 \subset r_m(\partial h_0) \subset r_m(h_0) \subset r_m( h_1) \subset r_m(h_2)
\subset \dots$$ Since $\mathscr{M}(x_0, y_0)$ is finite and since $y_0 \not \in r_m(h_0)$, there exists $j
> 0$ such that $y_0 \in r_m(h_j)$. Thus the wall $r_m(\partial h_i)$ separates $y_0$ from $y_i$ for all $i \geq j$.
Since $y_0$ and $y_i$ both lie on the wall $\mu$, it follows that $\partial h_i$ meets $\mu$ for all $i \geq j$.

This argument holds for any $\mu \in P \backslash \{m\}$. In particular, if we choose $\mu$ such that $m$ and
$\mu$ are separated by at least $L$ elements of $P$, where $L$ is the constant of Lemma~\ref{lem:triangle1}, we
deduce from this lemma that $W(\{m, \mu, \partial h_j\})$ is a Euclidean triangle subgroup. By
Corollary~\ref{cor:triangle2}, we obtain $r_{\partial h_i} \in  \overline{W(\{m, \mu, \partial h_j\})}$ for all
$i \geq j$. As before, this implies that $W(P') < \overline{W(\{m, \mu, \partial h_j\})}$ and, in particular,
that $m'' = r_m(m') = r_m(\partial h_0)$ meets $\mu$. Thus we have $\mu \in C$.

Since this holds for all walls $\mu \in P$ which are sufficiently far apart from $m$, and since $C$ is convex,
we finally deduce that $C = P$. By Lemma~\ref{lem:triangle1} this implies that $W(P \cup \{m'\})$ is a Euclidean
triangle subgroup. By Lemma~\ref{lem:Meucl}(i), we have $P \subset \Meucl$ whence $m \in \Meucl$.
\end{proof}

The main results of this section are the following two
propositions.

\begin{prop}\label{prop:groupe-Meucl}
The group $W(\Meucl)$ is isomorphic to a direct product of
finitely many irreducible affine Coxeter groups.
\end{prop}
\begin{proof}
We claim that for all $m, m' \in \Meucl$, either $P_F(m) = P_F(m')$ or the groups $W(P_F(m))$ and $W(P_F(m'))$
centralize each other or $W(P_F(m) \cup P_F(m'))$ is a Euclidean triangle subgroup.



We first deduce the desired result from the claim. We know that $W(\Meucl)$ is isomorphic to a Coxeter group.
Let $W(\Meucl) = W_1 \times \dots \times W_k$ be the decomposition of $W(\Meucl)$ in its direct components.
Hence $W_i$ is an irreducible Coxeter group for each $i = 1, \dots, k$. Let $M_i$ denote the set of walls
$m\in\Meucl$ such that $r_{m}\in W_i$. We note that $\Meucl=M_1\sqcup \dots\sqcup M_k$ and $W_i=W(M_i)$.

We must prove that $W_i$ is affine. We record the following easy observations which follow from the fact that
the $W_i$'s are the irreducible components of $W(\Meucl)$:
\begin{enumerate}
\item If $m \in \Meucl$ is a wall such that $r_{m} \in W_i$, then $W(P_F(m)) \leq W_i$.


\item If $m, m' \in \Meucl$ are two walls such that $P_F(m) \neq P_F(m')$  
 and that $r_{m}$ and $r_{m'}$ both
belong to $W_i$, then there exists a sequence of walls $m = m_0, m_1, \dots, m_\ell = m'$ such that for each
$j$, one has $m_j\in M_i$, $r_{m_j} \in W_i$ and $r_{m_j}$ does not commute with $r_{m_{j-1}}$ (a priori the
order of $r_{m_j} r_{m_{j-1}}$ might be infinite).


\end{enumerate}

We show that, in view of the claim above, these two observations imply that for any wall $m\in \Meucl$ such that $r_m \in W_i$, one has
$$W(P_F(m)) \leq W_i \leq \widetilde{W(P_F(m))},$$ where $\widetilde{W(P_F(m))}$ is the irreducible affine Coxeter
group provided by Lemma~\ref{lem:triangle2}.

By the first observation we just have to check that $W_i \leq \widetilde{W(P_F(m))}$. Since $W_i=W(M_i)$ it is enough to show that $r_{m'}\in \widetilde{W(P_F(m))}$ for any $m'\in M_i$. For such an $m'$ we have a sequence of walls $m = m_0, m_1, \dots, m_\ell = m'$ such that for each $j$,
one has ${m_j} \in M_i$ and $r_{m_j}$ does not commute with $r_{m_{j-1}}$. We are going to show by induction that for each $\mu\in P_F(m_i)$ we have $r_\mu\in \widetilde{W(P_F(m))}$, which implies in particular $r_{m'}\in \widetilde{W(P_F(m))}$.

This is clearly true for $i=0$. Assume this is true for $P_F(m_{i-1})$, with $i>0$. Either $m_i\in
P_F(m_{i-1})$, thus $P_F(m_i)=P_F(m_{i-1})$ and we have nothing to prove. Or,  by the initial claim,
$r_{m_i}r_{m_{i-1}}$ has finite order $>2$ and $W(P_F(m_{i-1})\cup P_F(m_{i}))$ is a Euclidean triangle
subgroup. Since $r_{m_i}$ and $r_{m_{i-1}}$ do not commute it follows that $W(P_F(m_{i-1})\cup \{m_i\})$ is a
Euclidean triangle subgroup.  Thus by Lemma~\ref{lem:triangle2} we have $r_{m_i}\in \widetilde{W(P_F(m))}$. In
fact the same argument applies to any wall $\mu \in P_F(m_i)$, which ends the proof.

The inclusion $W_i \leq \widetilde{W(P_F(m))},$ is now established.
In particular $W_i$ is an infinite reflection subgroup of an
irreducible affine Coxeter group; hence it must be itself an affine Coxeter group, as desired.

\smallskip
It remains to prove the claim. Let $m, m' \in \Meucl$.

Suppose that $P_F(m) \neq P_F(m')$. Then by Lemma~\ref{lem:propertymeucl} $m$ meets $m'$.

If there exists $m'' \in P_F(m) \cap P_F(m')$ then, by
Lemma~\ref{lem:Meucl}(i), we have $m'' \in \Meucl$ which implies
that the elements of $P_F(m'')$ are pairwise disjoint (see Lemma~\ref{lem:propertymeucl}). Since $m''
\in P_F(m) \cap P_F(m')$, we have $\{m, m' \} \subset P_F(m'')$ and,
hence, $m=m'$ because $m$ meets $m'$. This contradicts the fact
that $P_F(m) \neq P_F(m')$, thereby showing that $P_F(m) \cap P_F(m')$ is
empty. In other words, $m$ meets every element of $P_F(m')$ and $m'$
meets every element of $P_F(m)$.

For all $\mu \in P_F(m)$ we have $\mu \in \Meucl$ by
Lemma~\ref{lem:Meucl}(i) and, hence, $P_F(m) = P_F(\mu)$ by Lemma~\ref{lem:propertymeucl}. Similarly,
for all $\mu' \in P_F(m')$, we have $P_F(m')=P_F(\mu')$. Therefore, we
deduce from the previous paragraph that every element of $P_F(m)$
meets every element of $P_F(m')$.

Suppose moreover that $W(P_F(m))$ does not centralize $W(P_F(m'))$. Then there exist $p \in P_F(m)$ and $p' \in
P_F(m')$ such that $r_{p}$ and $r_{p'}$ do not commute. Let $p'':= r_p(p')$.


Suppose $p''$ meets only finitely elements of the line of walls $P_F(m)$. Then there is a segment of walls
$(p_-,p_1,p_2,\dots,p_n,p^+)$ inside $P_F(m)$ such that $\{p_1,p_2,\dots,p_n\}$ is the set of walls of $P_F(m)$
which meet $p''$, and $p''$ is disjoint from $p_-$ and $p_+$. We let $x_-,x_+$ denote points in $p'\cap
p_-,p'\cap p_+$ respectively. Since $p$ separates $p_-$ from $p_+$ and $p'\setminus p=p'\setminus p''$ we deduce
that $p''$ separates $x_-$ from $x_+$. Thus $p''$ separates  $p_-$ from $p_+$. In particular since $p_-$ and
$p_+$ meet $F$ we have $p''\in \MF$ and clearly $p''\in P_F(p_-)$. As we have already observed we have
$P_F(p_-)=P_F(m)=P_F(p)$. Thus $p''\in P_F(p)$, contradiction.

Thus in fact $p''$ meets infinitely many elements of $P_F(m)$. By Lemma~\ref{lem:triangle1}, this shows that
$W(P_F(m) \cup \{p'\})$ is a Euclidean triangle subgroup. Similarly $W(P_F(m')\cup \{p\})$ is a Euclidean triangle
subgroup. The order of the product $r_pr_{n'}$ is thus independent of the wall $n'$ chosen in the line of walls
$P_F(m')$. It follows that for each $n'\in P_F(m')$ the reflections $r_p$ and $r_{n'}$ do not commute. Then by
Lemma~\ref{lem:triangle1} the subgroup $W(P_F(m) \cup \{n'\})$ is also a Euclidean triangle subgroup. By
Corollary~\ref{cor:triangle2} we now deduce that $r_{n'}\in \overline{W(P_F(m) \cup \{p'\})}$. Thus $W(P_F(m'))
\subset \overline{W(P_F(m) \cup \{p'\})}$, and in particular the group $W(P_F(m) \cup P_F(m'))$ is a Euclidean
triangle subgroup, which proves the claim.
\end{proof}

\begin{cor}\label{cor:Meucl}
For all $m \in \Meucl$ and $\gamma \in W(\Meucl)$, if $\gamma . m \cap m = \emptyset$ then $\gamma.m \in
\Meucl$.
\end{cor}
\begin{proof}
By assumption, the group $\la r_m, r_{\gamma.m} \ra$ is an infinite dihedral group which is contained in
$W(\Meucl)$. Therefore, since $W(\Meucl)$ is an affine Coxeter group by Proposition~\ref{prop:groupe-Meucl}, the
group $W(P_F(m) \cup \{\gamma.m\})$ is an infinite dihedral group and, by Lemma~\ref{lem:propertymeucl}(iii), we
have $r_{\gamma.m} \in W(P_F(m))$. Since $P_F(m)$ is a convex line of walls, we deduce finally that $\gamma.m
\in P_F(m) \subset \Meucl$.
\end{proof}

\begin{prop}\label{prop:flat-walls}
One the following assertions holds:
\begin{enumerate}
\item[(i)] There exists an infinite subset $M \subset
 \mathscr{M}(F)$ which satisfies the following conditions:
\begin{itemize}
\item For all $m, m' \in M$, either $m \cap F = m' \cap F$
 or $m \cap F \cap m' = \emptyset$;

\item The groups $W(M)$ and $W(  \mathscr{M}(F) \backslash M)$ centralize each other.
\end{itemize}

\item[(ii)] The group $W( \mathscr{M}(F))$ is isomorphic to an affine Coxeter group.
\end{enumerate}
\end{prop}
\begin{proof}
Assume first that $ \Meucl=\MF$. Then by
Proposition~\ref{prop:groupe-Meucl} property (ii) holds.

Assume now there exists $m \in \MF \backslash \Meucl$. Let $M$ be the set of all those elements of $\MF$ which
do not belong to $\Meucl$ and which are $F$--parallel to $m$. By Lemma~\ref{lem:Meucl}(ii), we have $P_F(m)
\subset M$; in particular $M$ is infinite. Let $m' \in \MF \backslash M$. If $m'$ is not $F$--parallel to $m$
then $r_{m'}$ centralizes $W(M)$ by Lemma~\ref{lem:Meucl}(iii). If $m'$ is $F$--parallel to $m$, then $m' \in
\Meucl$ since $m' \not \in M$. In view of Lemma~\ref{lem:Meucl}(iv), this implies that $r_{m'}$ centralizes
$W(M)$. This shows that the groups $W(M)$ and $W(\MF \backslash M)$ centralize each other. Thus property (i)
holds.
\end{proof}

\section{From geometric flats to free abelian groups}\label{section:coxrank}

Let $X$ be a combinatorially convex subcomplex of the Davis
complex $|W|_0$, and $\Gamma$ be a subgroup of $W$ which
stabilizes $X$ and whose induced action on $X$ is cocompact. The
distance function on $|W|_0$ is denoted by $d$.

\begin{lem}\label{lem:cocompact}
Let $\rho \subset X$ be any unbounded subset through a given point $x$, and let $
\mathscr{M}(\rho):= \bigcup_{y,z \in \rho}  \mathscr{M}(y,z)$ be the set of walls which separate points of
$\rho$. There exists a constant $K$ (depending on $\rho$ and $\Gamma$) with the following property: given any
positive real number $r$, there exists a chamber $c$ at distance at most $K$ from $x$ and an element $\gamma \in
\Gamma \cap W( \mathscr{M}(\rho))$ such that $c$ and $\gamma.c$ both meet $\rho$, and that $d(c, \gamma.c)
> r$.
\end{lem}
\begin{proof}

Recall that a combinatorially convex subcomplex is a ($\mathrm{CAT}(0)$ convex) union of chambers.

Let $\mathcal C(\rho)$ denote the set of chambers  of $X$ meeting $\rho$: thus $\rho $ is covered by the
chambers of $\mathcal C(\rho)$. Recall that $\Gamma$ has finitely many orbits on the set of all chambers of $X$.
Since $\rho$ is unbounded, the set $\mathcal{C}(\rho)$ is infinite and it follows that there exists a chamber
$c\in C(\rho)$ such that $\Gamma.c\cap \mathcal C(\rho)$ is infinite.

We write $\Gamma.c\cap \mathcal C(\rho) = \{\gamma_0.c,\gamma_1.c,\dots,\gamma_i.c,\dots\}$ (with $\gamma_0=1$).
We pick a point $x_i$ in each intersection $\rho\cap \gamma_i.c$.  By Lemma~\ref{lem:basic-Coxeter}  there
exists $g_i\in W(\mathscr{M}(x_0,x_i))$ such that $g_ix_0$ and $x_i$ lie in a common chamber. Thus
${g_i}^{-1}\gamma_i$ is an element of $W$ sending $c$ to a chamber meeting $c$. There are finitely many such
elements.

Thus up to extracting a subsequence we may suppose that the sequence $({g_i}^{-1}\gamma_i)_{i\ge 1}$ is
constant. Then for each $i$ the element $\gamma'_{i}={\gamma_i}{\gamma_1}^{-1}$ belongs to $\Gamma\cap
W(\mathscr{M}(\rho))$. And also $\gamma'_{i}$ sends the chamber $\gamma_1c$ to the chamber $\gamma_ic$. The
Lemma follows because the set of chambers $(\gamma_ic)_{i\ge 1}$ is infinite.
\end{proof}

As before, let $ \mathscr{M}(F)$ denote the set of all walls which separate points of $F$. Theorem~\ref{thmA} of
the introduction is a straightforward consequence of the following:

\begin{thm}\label{thm:free-abelian}
Let $F$ be a geometric flat which is isometrically embedded in $X$; let $n$ denote its dimension. Then the
intersection $\Gamma \cap W( \mathscr{M}(F))$ contains a free abelian group of rank $n$.
\end{thm}
\begin{proof}
By Selberg's lemma, the group $\Gamma$ has a finite index subgroup
which is torsion free. Since $\Gamma$ is cocompact on $X$,  any
finite index subgroup of $\Gamma$ is cocompact as well, hence we
may assume without loss of generality that $\Gamma$ is torsion
free.

The proof works by induction on the dimension $n$ of the flat $F$.
We may assume
that $n>0$.

\medskip
Suppose first that $ \mathscr{M}(F)$ possesses a subset $M$ which satisfies the conditions~(i) of
Proposition~\ref{prop:flat-walls}. Let then $m$ be any element of $M$ and set $F':=F \cap m$. By Lemma~\ref{0},
$F'$ is a geometric flat of dimension $n-1$.

Let $\rho$ denote any geodesic ray of $F$ meeting transversally infinitely many walls of $M$. Let $x$ denote the
origin of $\rho$, and let $x_n$ denote the unique point of $\rho$ with $d(x,x_n)=n$. By
Lemma~\ref{lem:convexity} there exists a point $z_n\in X$ such that $\mathscr{M}(x,x_n)=\mathscr{M}(x,z_n)\sqcup
\mathscr{M}(z_n,x_n)$, with $\mathscr{M}(x,z_n) = \mathscr{M}(x, x_n) \cap M$. Observe that the cardinality of
$\mathscr{M}(x,z_n)$ tends to infinity with $n$, and thus $d(x,z_n)\to+\infty$. There is a subsequence
$(z_{n_k})_{k\ge 0}$ such that the geodesic segment $[x,z_{n_k}]\subset X$ converges to a geodesic ray
$\rho'\subset X$ (with origin $x$). Note that for every $y\in\rho'$ we have $\mathscr{M}(x,y)\subset
\mathscr{M}(x,z_{n_k})$ for $k$ large enough. In particular $\mathscr{M}(x,y)\subset M$. Thus
$\mathscr{M}(\rho')\subset M$.

We now apply Lemma~\ref{lem:cocompact} to the ray $\rho'$ for some (large) positive real number $r>0$. We then
get a nontrivial  element $\gamma\in \Gamma\cap W(M)$. Observe that $\gamma$ must be of infinite order since
$\Gamma$ is torsion free.

It follows from the definition of $M$ that $\gamma$ centralizes $W( \mathscr{M}(F'))$. Furthermore, since $W(
\mathscr{M}(F'))$ is isomorphic to a Coxeter group and since the center of any Coxeter group is a torsion group
(this is well known and is a straightforward consequence of \cite[Exercise 1, p.132]{Hum90}), the intersection
$W( \mathscr{M}(F')) \cap \la \gamma \ra$ is trivial. We deduce that the group generated by $W(
\mathscr{M}(F'))$ together with $\gamma$ is isomorphic to the direct product $W( \mathscr{M}(F')) \times \la
\gamma \ra$. The desired result follows  by induction.

\medskip
Suppose now that assertion~(ii) of Proposition~\ref{prop:flat-walls} holds. Let $\mu_1$ be any element of $
\mathscr{M}(F)$. Again by Lemma~\ref{0} the intersection $\mu_1\cap F$ is a geometric flat of dimension $n-1$.
Note that  any flat $\Phi$ of dimension $\ge 1$ is unbounded and thus has $\mathscr{M}(\Phi)\neq\emptyset$. Thus
for each $i = 2, \dots, n$ we may choose successively
$$\mu_i \in  \mathscr{M}(\big(\bigcap_{j=1}^{i-1} \mu_j \big) \cap
F).$$ In view of Lemma~\ref{0}, the set $\big(\bigcap_{i=1}^{n}
\mu_i \big) \cap F$ consists of a single point $x$ of $F$ and
for each $i \in \{1, \dots, n\}$, the set
$$\lambda_i:= \big(\bigcap_{j \in \{1, \dots, n\} \backslash\{i\}}  \mu_j \big) \cap
F$$ is a geodesic line of $F$.

We need the following auxiliary result:
\begin{lem}\label{lem:separable}
$\Gamma$ has a finite index subgroup $\Gamma'$ such that for any wall $m$ and any chamber $c$ meeting $m$, if
$\gamma\in\Gamma'$ sends $c$ to a chamber meeting $m$, then $\gamma m=m$.
\end{lem}
\begin{proof}
It is enough to prove the Lemma when $\Gamma=W$. Recall that the stabilizer of a wall $m$ is the centralizer of
the involution $r_m$. Since $W$ is residually finite the centralizer $Z(r_m)$ is a separable subgroup, that is
to say $Z(r_m)$ is an intersection of finite index subgroups. (In any residually finite group $W$ the
centralizer of any element $g$ is separable. Indeed, for $x \notin C= Z_W(g)$, we have $[x,g]\neq 1$, thus there
is a finite quotient $\phi:W\rightarrow \bar G$ such that $[\phi(x), \phi(g)]\neq 1$. Then $\phi (x)\notin
Z_{\bar{G}}(\phi(g)) $ and the finite index subgroup  $\phi\inv Z_{\bar{G}}(\phi(g))$ separates $x$ from $C$.)

We fix some wall $m$  and claim that there is a finite index subgroup $W_m\subset W$ such that for any chamber
$c$ meeting $m$, if $\gamma\in W_m$ sends $c$ to a chamber meeting $m$, then $\gamma m=m$. The lemma will follow
since we may assume that $W_m$ is normal, and there are only finitely many orbits of walls under $W$.

Let $B_m$ be the subset of $W$ consisting of all those elements $\gamma \in W$ such that there exists a chamber
$c$ such that $m$ and $\gamma. c$ both meet $m$. Note that $B_m$  is invariant by left- and right-multiplication
under $Z(r_m)$. In fact it is a finite union of double classes: $B_m=Z(r_m)\sqcup
Z(r_m)\gamma_1Z(r_m)\sqcup\dots\sqcup Z(r_m)\gamma_kZ(r_m)$, where $\gamma_1,\dots,\gamma_k$ do not belong to
$Z(r_m)$ (the finiteness follows from the fact that $Z(r_m)$ acts co-finitely on the set of chambers meeting
$m$, and from the local compactness of the Davis complex). The claim follows if we take for $W_m$ any finite
index subgroup of $W$ containing the separable subgroup $Z(r_m)$ but none of the elements
$\gamma_1,\dots,\gamma_k$.
\end{proof}
By Lemma~\ref{lem:separable} we may assume  that for any wall $m$ and any chamber $c$ meeting $m$, if
$\gamma\in\Gamma$ sends $c$ to a chamber meeting $m$, then $\gamma m=m$. Note that this implies in particular
that if $\gamma m$ intersects $m$, then  $\gamma m=m$.

Let $r$ be any positive real number. For each $i$ we choose one of the two rays contained in $\lambda_i$ with
origin $x$, and denote it by $\rho_i$. For each $i \in \{1, \dots, n\}$, Lemma~\ref{lem:cocompact} provides a
chamber $c_i$ at distance at most $K_i$ of $x$, and an element $\gamma_i(r) \in W(\mathscr{M}(\lambda_i))\cap
\Gamma$ ($\subset W(\mathscr{M}(F))\cap \Gamma$), such that $c_i \cap\rho_i$ and $\gamma_i(r).c_i \cap \rho_i$
are both nonempty, and that $d(c_i, \gamma_i(r).c_i) > r$. Here $c_i$ and $\gamma_i(r)$ depend on $r$, but $K_i$
depends only on $\rho_i$. Note that $\gamma_i(r)$ is of infinite order because $\Gamma$ is torsion free.

It immediately follows from the fact that $\rho_i \subset \mu_j$ that each $\gamma_i(r)$ preserves $\mu_j$
($j\neq i$).

Since $x \in \rho_i \cap \mu_i$ but $\rho_i \not \subset \mu_i$, it follows from
Lemma~\ref{lem:wall-neighborhood} that there is a constant $r_i$ such that, given any point $y$ of $\rho_i$, if
$y$ is at distance at least $r_i$ from $x$, then $y$ is at distance larger than $K_i + D$ from $\mu_i$, where
$D$ is the diameter of a chamber. Therefore, for each $r \geq r_i$, we have $d(x, \gamma_i(r).c_i) \geq d(c_i,
\gamma_i(r).c_i) > r$ and hence any point on $\gamma_i(r).c_i \cap \rho_i$ is at distance larger than $K_i + D$
from $\mu_i$. Thus $\gamma_i(r).c_i$ is at distance larger than $K_i$ from $\mu_i$. On the other hand
$\gamma_i(r).c_i$ is at distance at most $K_i$ from $\gamma_i(r)\mu_i$, from which it follows that
$\gamma_i(r)\mu_i\neq\mu_i$ for all $r \geq r_i$. By the above, this yields $\gamma_i(r)\mu_i \cap \mu_i =
\emptyset$ for all $r \geq r_i$.

Let $a_i$ be the half-space bounded by $\mu_i$ and containing $\rho_i$. We define an element $\gamma_i$ as
follows.

If $\gamma_i(r_i) a_i\subset a_i$ we set $\gamma_i=\gamma_i(r_i)$.

If not, then we choose $r > r_i$ as follows. Note that $\gamma_i(r) \mu_i \in \Meucl $ for all $r$ by
Corollary~\ref{cor:Meucl}. In particular $\gamma_i(r_i) \mu_i$ meets $\rho_i$, but $\rho_i \not \subset
\gamma_i(r_i) \mu_i$ because $x \in \rho_i \cap \mu_i$ and $\mu_i \cap \gamma_i(r_i) \mu_i = \emptyset$. Thus,
by Lemma~\ref{lem:wall-neighborhood}, every point of $\rho_i$ sufficiently far away from $x$ is also far way
from $\gamma_i(r_i)\mu_i$. Repeating the arguments used to define the constant $r_i$, we obtain a constant $r >
r_i$ such that $\gamma_i(r) \mu_i \neq \gamma_i(r_i) \mu_i$.

Now, if $\gamma_i(r) a_i \subset a_i$ we set  $\gamma_i=\gamma_i(r)$. Otherwise we set $\gamma_i=\gamma_i(r)\inv
{\gamma_i(r_i)}$. Let us check that, in the latter case, we have also $\gamma_i a_i\subset a_i$. The walls
$\mu_i$, $\gamma_i(r_i)\mu_i$ and $\gamma_i(r)\mu_i$ belong to $\Meucl$ by Corollary~\ref{cor:Meucl} and are
pairwise disjoint by construction. Thus they form a chain and it follows that $\gamma_i(r_i)a_i \subset
\gamma_i(r)a_i$ whence $\gamma_i a_i \subset a_i$. Therefore, for all $m > 0$, we have $\gamma_i ^m a_i\subset
a_i$ and hence $\gamma_i ^m \mu_i\cap \mu_i=\emptyset$ while $\gamma_i^m \mu_j = \mu_j$ for $j \neq i$.

Choose integers $m_1,\dots,m_n$ divisible enough so that each $\gamma'_i:={\gamma_i}^{m_i}$ belongs to the
translation subgroup of the affine Coxeter group $W( \mathscr{M}(F))$. Thus the $\gamma'_i$'s generate an
abelian group. In view of the action of each $\gamma'_i$ on the walls $\mu_1,\dots,\mu_n$, the intersection $\la
\gamma'_i \ra \cap \la \gamma'_j | \; j \neq i \ra$ is trivial for all $i$. This implies that the $\gamma'_i$'s
generate a free abelian group of rank~$n$.
\end{proof}

We note that the complete proof of Theorem
~\ref{thm:free-abelian} is much shorter when $(W,S)$ is assumed to be right-angled (in this case $\Meucl$ is empty).

\section{Geometric flats in Tits buildings}\label{section:buildrank}

The purpose of this section is to prove Theorem~\ref{thmC} of the
introduction.

As before, let $(W, S)$ be a Coxeter system of finite rank. Let $ \mathscr{B}=(\mathcal{C}( \mathscr{B}),
\delta)$ be a building of type $(W,S)$. Recall that $\mathcal{C}( \mathscr{B})$ is a set whose elements are
called \textbf{chambers}, and that $\delta: \mathcal{C}( \mathscr{B}) \times \mathcal{C}( \mathscr{B}) \to W$ is
a mapping, called \textbf{$W$-distance}, which satisfies the following conditions, where $x,y \in \mathcal{C}(
\mathscr{B}) $ and $w = \delta(x,y)$:

\begin{description}
\item[Bu1] $w = 1$ if and only if $x = y$;

\item[Bu2] if $z \in \mathcal{C}( \mathscr{B}) $ is such that $\delta(y,z) = s \in S$, then $\delta(x,z) = w$ or
$ws$, and if, furthermore, $l(ws) = l(w) + 1$, then $\delta(x,z) = ws$;

\item[Bu3] if $s \in S$, there exists $z \in \mathcal{C}( \mathscr{B}) $ such that $\delta(y,z) = s$ and
$\delta(x,z) = ws$.
\end{description}
For example the map $W\times W\to W$ sending $(x,y)$ to $x\inv y$ satisfies the above. An \textbf{apartment} of
the building $B$ is a subset $\mathcal{C}( \mathscr{A})\subset \mathcal{C}( \mathscr{B})$ such that there exists
a bijection $f:\mathcal{C}( \mathscr{A})\to W$ satisfying $\delta(x,y)=f(x)\inv{f(y)}$.

The composed map $\ell \circ \delta : \mathcal{C}( \mathscr{B}) \times \mathcal{C}( \mathscr{B}) \to \N$, where
$\ell$ is the word metric on $W$ with respect to $S$, is called the \textbf{numerical distance} of
$\mathscr{B}$. It is a discrete metric on $\mathcal{C}( \mathscr{B})$.

The following lemma is well known:

\begin{lem}\label{lem:apt:extension}
Let $\mathcal{C}( \mathscr{A})$ be an apartment and $C$ be a subset of $\mathcal{C}( \mathscr{B})$. Suppose that
there exists a map $f: C \to \mathcal{C}( \mathscr{A})$ such that $\delta(f(c), f(d))=\delta(c, d)$ for all $c,
d \in C$. Then there exists an apartment $ \mathcal{C}( \mathscr{A}')$ such that $C \subset \mathcal{C}(
\mathscr{A}')$.
\end{lem}
\begin{proof}
Follows from \cite[\S 3.7.4]{Ti81}.
\end{proof}

Let $T\subset S$ and let $c$ be a chamber of the building $B$. The \textbf{residue of type $T$ of $c$} is the
set $\rho_T(c)$ of those chambers $c'$ for which $\delta(c,c')\in W(T)$. The residue is called
\textbf{spherical} whenever $W(T)$ is finite. Given any residue $\rho$ of $B$ and any chamber $x$, there exists
a unique chamber $c$ in $\rho$ at minimal numerical distance from $x$. This chamber has the property that
$\delta(x,d)=\delta(x,c)\delta(c,d)$ for each chamber $d$ of $\rho$. The chamber $c$ is called the
\textbf{projection of $x$ onto $\rho$} and is denoted by $\proj_\rho(c)$ (see \cite[\S Corollary~3.9]{Ro89}).

\begin{lem}\label{lem:apt:criterion}
Let $\mathcal{C}( \mathscr{A})$ be an apartment of $ \mathscr{B}$ and $C \subset \mathcal{C}( \mathscr{A})$ be a
set of chambers. Suppose that there exists a residue $\rho$ and a chamber $c \in C$ such that $c \in
\mathcal{C}(\rho)$ and $\proj_\rho(c')=c$ for all $c' \in C$. Then, for any chamber $d \in \mathcal{C}(\rho)
\backslash \{c\}$, there exists an apartment $ \mathcal{C}( \mathscr{A}_d)$ such that $C \cup \{d\}$ is
contained in $\mathcal{C}( \mathscr{A}_d)$.
\end{lem}
\begin{proof}
Let  $d \in \mathcal{C}(\rho) \backslash \{c\}$ and let $w_d:= \delta(c,d)$. Let $d'$ be the unique chamber of
$\mathcal{C}( \mathscr{A})$ such that $\delta(c,d')=w_d$. For any $c' \in C$, we have
$\delta(c',d)=\delta(c',c).w_d = \delta(c',d')$ because $\proj_\rho(c')=c$. It follows that the function $f: C
\cup \{d\} \to  C \cup \{d'\}$, which maps $d$ to $d'$ and induces the identity on $C$, preserves the
$W$-distance $\delta$. Therefore, the existence of an apartment $\mathcal{C}(  \mathscr{A}_d)$ such that
$\mathcal{C}( \mathscr{A}_d)$ contains $C \cup \{d\}$ follows from Lemma~\ref{lem:apt:extension}.
\end{proof}

Before stating the main result of this section, we need to introduce some additional terminology and notation:
\begin{itemize}
\item $| \mathscr{B}|_0$ denotes the $\mathrm{CAT}(0)$-realization of the building $ \mathscr{B}$, as defined in
\cite{Da98}; it is a piecewise Euclidean simplicial complex. For each chamber $c\in B$ there is an associated
$\mathrm{CAT}(0)$-convex subcomplex $| {c}|_0 \subset | \mathscr{B}|_0$, which we call  \textbf{the associated
geometric chamber}. For every subset $C\subset \mathcal{C}(\mathscr{B})$ we denote by $| {C}|_0$ the union of
geometric chambers $| {c}|_0$ associated to chambers $c\in C$. We say that a subcomplex $X\subset |
\mathscr{B}|_0$ is \textbf{combinatorial} whenever it is a union of geometric chambers. If $\mathscr{A}$ is any
apartment of $\mathscr{B}$ the subcomplex $| \mathscr{A}|_0$ is isometric to $ | {W}|_0$. As a simplicial
complex, $| \mathscr{A}|_0$ is isomorphic to the first barycentric subdivision of the Davis complex $ | {W}|_0$.

\item Given $x \in | \mathscr{B}|_0$, we set
$$\rho(x) := \{c \in \mathcal{C}(\mathscr{B}) | \; x \in |c|_0\}$$
and $$\sigma(x) := \bigcap_{c \in \rho(x)} |c|_0.$$ The set $\rho(x)$ is a (spherical) residue. The subcomplex
$| {\rho(x)}|_0$ is a neighbourhood $N(x)$ of $x$ in $ | \mathscr{B}|_0$. For every chamber $c \in
\mathcal{C}(\mathscr{B})$, the set ${\rm Int}(c)$ of points $x\in | \mathscr{B}|_0$ such that $\rho(x)=\{c\}$ is
an open subset of $| \mathscr{B}|_0$. It is the interior of $| {c}|_0$ and its closure is $| {c}|_0$.

\item Given $E \subset | \mathscr{B}|_0$, we set
$$\mathcal{C}(E) := \{c \in \mathcal{C}(\mathscr{B}) | \; |c|_0 \subset E\}.$$
For example given any $x\in  | \mathscr{B}|_0$ we have $\mathcal{C}(N(x))=\rho(x)$. We say that a subcomplex
$\mathscr{A} \subset  | \mathscr{B}|_0$ is a \textbf{geometric apartment} provided $\mathscr{A}$ is
combinatorial and $\mathcal{C}(\mathscr{A})$ is an apartment of $\mathscr{B}$.

\item Given a geometric flat $F \subset | \mathscr{B}|_0$ and any subset $E \subset | \mathscr{B}|_0$, we denote
by $\dim(F\cap E)$ the dimension of the Euclidean subspace of $F$ generated by $E \cap F$; by convention, the
empty set is a Euclidean subspace of dimension $-1$.
\end{itemize}

Let now $F\subset | \mathscr{B}|_0$ be a geometric flat of dimension $n$. Since the combinatorial subcomplexes
$N(x)$ are neighborhoods of $x$, we have:
$$
\forall x \in F, \ \exists\ c \in \mathcal{C}( \mathscr{B}) \text{ such that } x \in |c|_0 \text{ and }
\dim(F\cap |c|_0) = n.
$$
And since every geometric chamber is the closure of its interior, we deduce:
$$
\forall x \in F, \ \exists\ y \in F \text{ such that } x \in
\sigma(y) \text{ and } \dim(F\cap \sigma(y))= n.
$$
These two basic facts will be used repeatedly in the following.

\begin{thm}\label{thm:flats-in-apts}
Let $F\subset | \mathscr{B}|_0$ be a geometric flat of dimension $n$ and let $c_0$ be a chamber such that
$\dim(F\cap c_0)=n$ (the  geometric chamber associated to $c_0$ is also denoted by $c_0$). Define
$$C(F,c_0):= \{\proj_{\rho(x)}(c_0) \; | \; x \in F\}.$$
Then there exists a geometric apartment $\mathscr{A}$ such that $C(F,c_0) \subset \mathcal{C}( \mathscr{A})$. In
particular, we have $F \subset  \mathscr{A}$.
\end{thm}
\begin{proof}
The proof is by induction on $n$, the case $n=0$ being trivial. We
assume now that $n>0$.

Let $F_0 \subset F$ be a Euclidean hyperplane such that $\dim({F_0} \cap c_0) = n-1$. By induction, the set
$C(F_0, c_0)$ is contained in the set of chambers of some apartment. In view of Lemma~\ref{lem:apt:extension},
it follows from Zorn's lemma that the collection of all those subsets of $C(F, c_0)$ which contain $C(F_0, c_0)$
and which are contained in the  set of chambers of some apartment, has a maximal element.

Let $C_1$ be such a maximal element and choose a geometric apartment $\mathscr{A}_1$ such that $C_1 \subset
\mathcal{C}( \mathscr{A}_1)$. Set $X:=  \mathscr{A}_1 \cap F$. Note that $X$ is closed and convex.

Suppose by contradiction that $C_1$ is properly contained in $C(F,
c_0)$. The rest of the proof is divided into several steps. The
final claim below contradicts the maximality of $C_1$, thereby
proving the theorem.

\setcounter{claim}{0}
\begin{claim}\label{claim1}
For all $x \in X$, we have $\proj_{\rho(x)}(c_0) \in C_1$.
\end{claim}

Since $ \mathscr{A}_1$ is a combinatorial subcomplex, we have $\sigma(x) \subset
 \mathscr{A}_1$. Since $c_0 \in \mathcal{C}( \mathscr{A}_1)$, we
have $\proj_{\rho(x)}(c_0) \in \mathcal{C}( \mathscr{A}_1)$. Therefore, the claim follows from the maximality of
$C_1$.

\begin{claim}\label{claim2}
For all $c \in C_1$, there exists $x \in X$ such that
$\proj_{\rho(x)}(c_0)=c$.
\end{claim}

Given $c \in C(F,c_0)$, there exists $x \in F$ such that $\proj_{\rho(x)}(c_0)=c$. If now $c \in C_1$, then
$\sigma(x) \subset | {c}|_0 \subset  \mathscr{A}_1$. Thus $x \in F \cap
 \mathscr{A}_1=X$.

\begin{claim}\label{claim3}
$\dim(F\cap X) = n$.
\end{claim}

Clear since $c_0 \cap F \subset X$ and $\dim(F\cap c_0)=n$.

\begin{claim}\label{claim4}
There exists a Euclidean hyperplane $F_1 \subset F$ which is contained in $ \mathscr{A}_1$ and which bounds an
open half-space of $F$, none of whose points is contained in $ \mathscr{A}_1$. In other words, the hyperplance
$F_1$ is contained in the Euclidean boundary $\partial X$ of $X$.
\end{claim}

Let $c \in  C(F, c_0) \backslash C_1$ and let $x \in F$ be such that $\proj_{\rho(x)}(c_0)=c$. By
Claim~\ref{claim1}, $x$ does not belong to $X$. Given $x_0 \in c_0 \cap F$, we have $[x_0,x] \cap X = [x_0, y]$
for some $y \in X$ because $X$ is closed and convex. Let $F_1 \subset F$ be the Euclidean hyperplane parallel to
$F_0$ and containing $y$. We have $F_1 \subset X$ by convexity. Furthermore, it is clear from the definition of
$y$ and $F_1$ that any point $z \in F\backslash F_1$ on the same side of $F_1$ as $x$ does not belong to $
\mathscr{A}_1$.

\begin{claim}\label{claim5}
Let $x_1 \in F_1$ be such that $\dim({F_1}\cap \sigma(x_1))=n-1$. For
all $c \in C_1$, we have $\proj_{\rho(x_1)}(c)=
\proj_{\rho(x_1)}(c_0)$.
\end{claim}

Let $c_1:= \proj_{\rho(x_1)}(c_0)$. Suppose by contradiction that there exists $c \in C_1$ such that
$\proj_{\rho(x_1)}(c) \neq c_1$. Let $h$ be a (Coxeter) half-space of the apartment $
\mathcal{C}(\mathscr{A}_1)$ containing $c_1$ but not $c_2:=\proj_{\rho(x_1)}(c)$. Thus $h$ contains $c_0$ but
not $c$.

Since $\sigma(x_1) \subset | {c_1}|_0 \cap  | {c_2}|_0$, we have $\sigma(x_1)
\subset
\partial  | {h}|_0$. Therefore, since $F_1 \subset
\mathscr{A}_1$ (see Claim~\ref{claim4}) and since $\dim({F_1}\cap \sigma(x_1))=n-1$, we deduce from
Lemma~\ref{0} that $F_1 \subset \partial  | {h}|_0$. By Claim~\ref{claim4}, the set $X$, as a subset of $F$, is
entirely contained in one of the Euclidean half-spaces of $F$ determined by $F_1$. Since $F_1 \subset \partial |
{h}|_0$, we deduce that $X$, as a subset of $\mathscr{A}_1$, is entirely contained in one of the Coxeter
half-spaces of $\mathscr{A}_1$ determined by $ \partial  | {h}|_0$. Since $c_0 \subset X \cap |h|_0$, we obtain
$X \subset   | {h}|_0$.

Since $c \in C_1$, there exists $x \in X$ such that $\proj_{\rho(x)}(c_0)=c$ by Claim~\ref{claim2}. Since $X
\subset  | {h}|_0$ and since $ | {h}|_0$ is a combinatorial subcomplex, we have $\sigma(x) \subset  | {h}|_0$
and hence $\proj_{\rho(x)}(c_0) \in h$ by the combinatorial convexity of Coxeter half-spaces. This contradicts
the fact that $h$ does not contain $c$.

\begin{claim}
There exists $d \in C(F, c_0)$ and an apartment $ \mathscr{A}_d$ such that $C_1 \cup \{d\} \subset \mathcal{C}(
\mathscr{A}_d)$.
\end{claim}

Let $x_1 \in F_1$ be as in Claim~\ref{claim5}. By Claim~\ref{claim1} we have $c_1:=\proj_{\rho(x_1)}(c_0)\in
C_1$.  Let $y \in F \backslash X$ be such that $x_1 \in \sigma(y)$. Let $d:= \proj_{\sigma(y)}(c_0)$. Clearly $d
\in C(F, c_0)$. Furthermore $d \not \in C_1$, otherwise we would have $y \in \sigma(y) \subset d \subset
\mathscr{A}_1$, whence $y \in X$, which is absurd. Since $\sigma(x_1) \subset \sigma(y) \subset d$, the claim
follows from Lemma~\ref{lem:apt:criterion} together with Claim~\ref{claim5}.
\end{proof}

Clearly, Theorem~\ref{thmC} of the introduction is an immediate
consequence of Theorem~\ref{thm:flats-in-apts}, combined with
Corollary~\ref{corB}.


\vspace{1cm}

\addtolength{\parindent}{-1.6pt}
\end{document}